\documentclass[10pt,twoside]{siamltex}
\usepackage{amsfonts,epsfig}

\newtheorem{remark}[theorem]{Remark}

\newtheorem{question}[theorem]{Question}
\newtheorem{prop}[theorem]{Proposition}
\newtheorem{cor}[theorem]{Corollary}
\newcommand\set[1]{\{#1\}}

\begin{document}



\bibliographystyle{plain}
\title{Clifford algebras of $p$-central sets}

\author{Adam\ Chapman\thanks{Department of Mathematics,
Bar Ilan University, Ramat-Gan, Israel
(adam1chapman@yahoo.com). Ph.D student under the supervision of Prof. Uzi Vishne}}


\pagestyle{myheadings}
\maketitle

\begin{abstract}
A generalization of the term ``generalized Clifford algebras" (as appears in papers on advances in applied Clifford algebras) is introduced. This algebra is studied by means of structure theory of central simple algebras.
A graph theoretical approach is proposed for studying the generating set of this algebra in case where the prime number under discussion is three.
Finally, it is shown how to obtain solutions in to the equation $\alpha Y^3=\alpha X_1^3+\beta X_2^3+\alpha^2 \beta^2 X_3^3$ in $\mathbb{Z}[\rho]$ where $\rho$ is the primitive third root of unity.
\end{abstract}


\section{Introduction} \label{intro-sec}
\subsection{The Clifford algebra of a $p$-central set}
Let $F$ be a field containing a primitive $p$th root of unity $\rho$, and $A$ be an associative $F$-algebra.
A finite subset $B={b_1,\dots,b_n} \subseteq A$ consisting of invertible elements is called a \textbf{$p$-central set}\footnote{This term was introduced by Rowen in \cite[Vol II, pp. 248-251]{Rowen}} if
\begin{enumerate}
\item For any $1 \leq k \leq n$, $b_k^p=\alpha_k \in F$.
\item For any $1 \leq m < k \leq n$, $b_m b_k=\rho^{d_{m,k}} b_k b_m$
\end{enumerate}

The Clifford algebra of $B$, denoted by $C(B)$, is defined to be $F[x_1,\dots,x_n : x_k^p=\alpha_k^p, x_m x_k=\rho^{c_{m,k}} x_k x_m$.
This generalizes the definition of a generalized Clifford algebra as appears in papers on advances in applied Clifford algebras, where all $d_{m,k}$-s are equal to one. (See for example \cite{Koc})

The subalgebra $F[b_1,\dots,b_n]$ is obviously a homomorphic image of $C(B)$.
Since the Clifford algebra is depended only on the choice of $\alpha$-s and $c$-s, we can define the Clifford algebra according to these elements instead of referring to a specific $p$-central set $B$.

In \cite[Vol. II, pp. 248-251]{Rowen} it is explained that $F[b_1,\dots,b_n]$ decomposes as a tensor product of $p$-cyclic algebras over some extension of $F$.
Furthermore, a way of how to obtain this decomposition is suggested:
One should change the $p$-generating set by replacing some $b_i$ with $b_i b_j^k$ (for some $1 \leq k \leq p-1$) until he receives a set of ``best possible form", i.e. a set whose center is the biggest possible, say $m$. Then this algebra decomposes as a tensor product of $\frac{n-m}{2}$ $p$-cyclic algebras and the field extension generated by the center of the generating set of best possible form.

It does not say, however, how can one obtain this $p$-central set of best possible form.

In Section \ref{structure}, a detailed arithmetical method is given for completely analyzing the structure of $C(B)$ (and consequently the structure of $F[b_1,\dots,b_n]$).

\subsection{Coherent $p$-sets}

A subspace $V=F x_1+\dots+F x_n \subset A$ is called \textbf{$p$-central} if each $v \in V$ satisfies $v^p \in F$.
If its basis $\set{x_1,\dots,x_n}$ forms a $p$-central set then we call it a \textbf{coherent $p$-central set}.

\begin{remark}
It should be mentioned that if $\set{x_1,\dots,x_n}$ is a coherent $p$-central set and $x_i$ commutes with $x_j$ then $F x_i=F x_j$. Consequently, if this set spans a $p$-central space of dimension $n$ then it contains no commuting pair of elements.
\end{remark}

\begin{question}
When is a $p$-central set coherent?
\end{question}

For $p=2$ it is always true, i.e. every $2$-central set is coherent. (See \cite{Lam} for further details)

In order to answer this question we turn to graph theory.
A directed graph is a pair of sets $(V,E)$. $V$ is the set of vertices and $E$ is the set of edges.
Every edge is an ordered pair $(a,b)$ where $a,b \in V$.
A cycle is a sequence of vertices $a_1,a_2,\dots,a_n$ with no repetitions such that $(a_1,a_2),\dots,(a_{n-1},a_n),(a_n,a_1) \in V$.

Later in this paper will be shown how coherent $3$-central sets have unique graphical form.

\subsection{Solving a Diophantine equation}

In \cite{Aragon} it is shown how Clifford algebras can be used for solving quadratic Diophantine equations such as $\alpha Y^2=\alpha X_1^2+\beta X_2^2$.

Here we show how a similar technique provides solutions to the cubic Diophantine equation $\alpha Y^3=\alpha X_1^3+\beta X_2^3+\alpha^2 \beta^2 X_3^3$ in $\mathbb{Z}[\rho]$ (for $p=3$).

For the case of $\alpha=\beta=1$, these solutions are completely distinct from the known integral solutions, such as Euler's (See \cite[Vol. II, pp. 550-561]{Dickson}).

\section{The structure of the Clifford algebra of a $p$-central set}\label{structure}
Let $C=F[x_1,\dots,x_n : x_k^p=\alpha_k^p, x_m x_k=\rho^{c_{m,k}} x_k x_m]$.
In particular, $B=\set{x_1,\dots,x_n}$ is a $p$-central set and it generates $C$ therefore we shall call it a generating $p$-central set.
We build a matrix as follows:
$M_B=\left[ \begin{array}{ccc}
c_{1,1} & \ldots &  c_{1,n}  \\
\vdots & \ldots & \vdots \\
c_{n,1} & \ldots & c_{n,n}  \end{array} \right] \in M_n(\mathbb{Z}/p \mathbb{Z})$
and a vector $v=[\alpha_1,\dots,\alpha_n]$.
As a matrix in $M_n(\mathbb{Z}/p\mathbb{Z})$, $M_B$ is skew-symmetric, and therefore similar to some block matrix $B=(\oplus_{k=1}^m H) \oplus 0_{(n-2 m) \times (n-2 m)}$ where $H=\left[\begin{array}{lr}
0 & -1\\
1 & 0 \end{array}\right]$.
Consequently there exists an orthogonal matrix $D \in M_n(\mathbb{Z}/p\mathbb{Z})$ for which
$M'=D M D^t$.

Without changing the algebra $C$, we can change its generating $p$-central set $B=\set{x_1,\dots,x_n}$ by replacing each $x_i$ with $x_1^{d_{i,1}} x_2^{d_{i,2}} \dots x_n^{d_{i,n}}$ and obtain an altered generating $p$-central set $B'$.
This change affects the matrix $M_B$ who is in turn replaced with $M_{B'}=D M_B D^t$, where $D=(d_{i,j})$.

On the other hand, any change of matrices $M_B \rightarrow D M_B D^t$ can be obtained by a change of generating $p$-central sets $\set{x_1,\dots,x_n} \rightarrow \set{x_1^{d_{1,1}} \dots x_n^{d_{1,n}},\dots,x_1^{d_{n,1}} \dots x_n^{d_{n,n}}}$.

In particular, because there exists an orthogonal matrix $D$ such that $M'=D M_B D^t$, $M'=M_{B'}$ is the matrix related to the generating $p$-central set $B'=\set{x_1^{d_{1,1}} \dots x_n^{d_{1,n}},\dots,x_1^{d_{n,1}} \dots x_n^{d_{n,n}}}$.

Denoting $y_i=x_1^{d_{i,1}} x_2^{d_{i,2}} \dots x_n^{d_{i,n}}$ for each $1 \leq i \leq n$, the algebra $C$ is isomorphic to $F[y_1,y_2] \otimes \dots \otimes F[y_{2 m-1},y_{2 m}] \otimes F[y_{2 m+1},\dots,y_n]$.
For each $1 \leq k \leq m$, $F[y_{2 k-1},y_{2 k}]$ is the $p$-cyclic algebra $(\alpha_1^{d_{2 k-1,1}} \dots \alpha_n^{d_{2 k-1,n}},\alpha_1^{d_{2 k,1}} \dots \alpha_n^{d_{2 k,n}})_{p,F}$, and $F[y_{2 m+1},\dots,y_n : y_i^p=\alpha_1^{d_{i,1}} \dots \alpha_n^{d_{i,n}} \forall 2 m+1 \leq i \leq n]$ is the commutative ring (maybe a field, but not necessarily).

The algebra $C$ is obviously Azumaya and each simple homomorphic image of it is of degree $p^m$ and exponent either $p$ or $1$.

\section{Coherent $3$-central sets}
Let $A$ be a central simple algebra over the field $F$ containing a primitive 3rd root of unity $\rho$. All the discussion below takes place in this algebra.

According to \cite[Corollary 2.2]{Chapman}, a set $\set{x_1,\dots,x_n}$ spans a 3-central spaces if and only if every subset of cardinality three $\set{x_i,x_j,x_k}$ spans a 3-central space.
Therefore we will start with the set of cardinality 3.

Let $\set{x,y,z}$ span a 3-central space such that $x \ast y \ast z=0$.
Now, $y=y_1+y_2$ and $z=z_1+z_2$ such that $y_i x=\rho^i x y_i$ and $z_i x=\rho^i x z_i$.
Consequently
$x \ast y \ast z=x ((1-\rho^2)(y_1 z_2+z_1 y_2)+(1-\rho) (y_2 z_1+z_2 y_1))=0$, which means that
$y_1 z_2+z_1 y_2=\rho (y_2 z_1+z_2 y_1)$.

If $\set{x,y,z}$ is also a 3-central set then either $y_1=0$ or $y_2=0$ and either $z_1=0$ or $z_2=0$. Henceforth (up to change of order of the elements) $y x=\rho x y$, $z y=\rho y z$ and $z x=\rho x z$.

\begin{question}
What happens if we do not require $x \ast y \ast z=0$?
\end{question}

For a 3-central set $\set{x,y,z}$ spanning a 3-central space there are two options:

\begin{enumerate}
\item \label{fc} $x \ast y \ast z=0$ which is the case we already dealt with.
\item \label{sc} $x \ast y \ast z \neq 0$, but still $x \ast y \ast z \in F$. Because $\set{x,y,z}$ is a 3-central set, $x \ast y \ast z=a x y z$ where $a \in F$. Consequently $z \in F x^{-1} y^{-1}$. Since multiplying an element by a central element does not change any of the characteristic we are dealing with here, we shall say that in this case $z=x^{-1} y^{-1}$.
\end{enumerate}

Before we proceed, there are two relevant matrices to each 3-central set:
Given a 3-central set $B=\set{x_1,\dots,x_n}$, there is the matrix $M_B$ as in Section \ref{structure}.
This matrix has an underlying directed graph $(B,E_B)$ such that $(x_i,x_j) \in E \Leftrightarrow a_{i,j}=1$.
Disregarding the orientation, this graph is complete.

Every cycle of length 3 in this graph is of Case (\ref{sc}), i.e. three elements whose product is central.

\begin{prop}
Let $(x_i,x_j,x_k)$ be a cycle, then for every $x_l$ outside this cycle, either $(x_l,x_i),(x_l,x_j),(x_l,x_k) \in E_B$ or $(x_i,x_l),(x_j,x_l),(x_k,x_l) \in E_B$.
\end{prop}

\begin{proof}
If $(x_l,x_i),(x_j,x_l) \in E_B$ then $(x_l,x_i,x_j)$ is a cycle, and so $x_l=x_k$ (up to multiplication by a central element), which is not true. Consequently if $(x_l,x_i) \in E_B$ then also $(x_l,x_j) \in E_B$ and inductively also $(x_l,x_k) \in E_B$.

Similarly, if $(x_j,x_l) \in E_B$ then also $(x_i,x_l),(x_k,x_l) \in E_B$.
\end{proof}

\begin{prop}
No cycle shares a vertex with another.
\end{prop}

\begin{proof}
Let $(x_i,x_j,x_k)$ and $(x_i,x_r,x_s)$ be two cycles. If $x=r$ then $k=s$ which makes them the same.
So, we assume that $j \neq r$ and $k \neq s$.
If $(x_j,x_s) \in E_B$ then $(x_i,x_j,x_s)$ is a cycle, and then $s=k$, i.e. a contradiction.
Similarly we will have a contradiction if $(x_s,x_j) \in E_B$.
\end{proof}

\begin{prop}
There are no cycles of length greater than 3.
\end{prop}

\begin{proof}
Let $x_1,\dots,x_r$ be a cycle of length $r \geq 4$.
Let $i$ be the maximal index for which $(x_1,x_i) \in E_B$.
It exists because $(x_1,x_2) \in E_B$ and it is no more than $r-1$ because $(x_r,x_1) \in E_B$.
Now, $(x_{i+1},x_1) \in E_B$.
Therefore $(x_1,x_i,x_{i+1})$ is a cycle.
If $i \geq 3$ then $x_{i-1}$ is outside this cycle, therefore $(x_{i-1},x_1) \in E_B$. Let $j$ be the minimal index for which $(x_{j+1},x_1) \in E_B$. $j \neq i$, and still $(x_1,x_j) \in E_B$.
Therefore we have another cycle $(x_1,x_j,x_{j+1})$ which shares a vertex with the previous one, and that is impossible.
If $i=2$ then $x_4$ is outside of the circle and since $(x_3,x_4) \in E_B$, we have $(x_1,x_4) \in E_B$ which contradicts the maximality of $i$.
\end{proof}

Given $(B,E_B)$ we can diminish it as follows:
Of each cycle $(x_i,x_j,x_k)$, we omit one of the vertices arbitrarily from $B$ and its related edges from $E_B$.
Finally we are left with a graph $(\underline{B},E_{\underline{B}})$ with no cycles.

Now, the algebra generated by $\underline{B}$ is isomorphic to a tensor product of $\lfloor \frac{\sharp\underline{B}}{2} \rfloor$,
and since the algebra generated by $B$ is the same as the algebra generated by $\underline{B}$, we obtain the following easy result:

\begin{cor}
The maximal coherent 3-central set in a tensor product of $m$ cyclic algebras of degree 3 is of cardinality $3 m+1$.
\end{cor}

\section{Solving a Diophantine equation}

Let $A$ be a central simple $\mathbb{Q}[\rho]$-algebra containing two $3$-central elements $x,y \subseteq $ satisfying $y x=\rho x y$, $x^3=\alpha$ and $y^3=\beta$, where $\alpha,\beta \in \mathbb{Z}[\sqrt{3},i]$.

$(a x+b y) x (a x+b y)^{-1}
=(a^3 \alpha+b^3 \beta)^{-1} (a x+b y) x (a^2 x^2-\rho^{-1} a b x y+b^2 y^2)
=(a^3 \alpha+b^3 \beta)^{-1} ((a^3 \alpha+\rho b^3 \beta) x+(1-\rho^{-1}) b a^2 \alpha y+(1-\rho) a b^2 x^2 y^2)$

The set $\set{x,y,x^2 y^2}$ is a coherent $3$-central set. In particular $x * y * (x^2 y^2)=x y x^2 y^2+x^3 y^3+y x^3 y^2+y x^2 y^2 x+x^2 y^2 x y+x^2 y^3 x=-3 \rho \alpha \beta$.
Consequently, $\alpha=((a x+b y) x (a x+b y)^{-1})^3=(a^3 \alpha+b^3 \beta)^{-3} ((a^3 \alpha+\rho b^3 \beta)^3 \alpha+(1-\rho^{-1})^3 a^6 b^3 \alpha^3 \beta+(1-\rho)^3 a^3 b^6 \alpha^2 \beta^2-9 \rho a^3 b^3 \alpha (a^3 \alpha+\rho b^3 \beta) \alpha \beta)$, which means that
$\alpha (a^3 \alpha+b^3 \beta)^3=\alpha (a^3 \alpha+\rho b^3 \beta)^3+3 (1-\rho) \beta (a^2 b \alpha)^3+3 (1-\rho^{-1})\alpha^2 \beta^2 (a b^2)^3$

This means that for any $a,b \in O$, $Y=a^3 \alpha+b^3 \beta$, $X_1=a^3 \alpha+\rho b^3 \beta$, $X_2=a^2 b \alpha$ and $X_3=a b^2$ are solutions in $\mathbb{Z}[\rho]$ (the ring of integers of the field $\mathbb{Q}[\rho]$) to the equation
$\alpha Y^3=\alpha X_1^3+3 (1-\rho) \beta X_2^2+3 (1-\rho^{-1}) \alpha^2 \beta^2 X_3^2$.

Since $\alpha$ can be replaced with $3 (1-\rho) \gamma$ we actually have solutions in $\mathbb{Z}[\rho]$ to any equation of the form
$\gamma Y^3=\gamma X_1^3+\beta X_2^3+27 \gamma^2 \beta^2 X_3^3$.
Now, $27=3^3$, and so we have solutions to the equation
$\gamma Y^3=\gamma X_1^3+\beta X_2^3+\gamma^2 \beta^2 X_3^3$ for any $\gamma$ and $\beta$.
In this case the solution will be $Y=3 (1-\rho) a^3 \gamma+b^3 \beta$, $X_1=3 (1-\rho) a^3 \gamma+\rho b^3 \beta$, $X_2=3 (1-\rho) a^2 b \gamma$ and $X_3=3 a b^2$.

By replacing $b$ with $(1-\rho) c$ our solutions to
$\gamma Y^3=\gamma X_1^3+\beta X_2^3+\gamma^2 \beta^2 X_3^3$ become $Y=a^3 \gamma-\rho \beta c^3$, $X_1^3=a^3 \gamma-\rho^2 \beta c^3$, $X_2=\gamma (1-\rho) a^2 c$ and $X_3=(1-\rho) a c^2$.


\end{document}